\documentclass[11pt]{article}

\usepackage{latexsym}
\usepackage{graphicx}
\usepackage{framed}
\usepackage{threeparttable}
\usepackage{amssymb}
\usepackage{amsmath}
\usepackage{multirow}
\usepackage{hyperref}

\begin{document}

\newcommand{\noi}{\noindent}
\newcommand{\nn}{\nonumber}
\newcommand{\bd}{\begin{displaymath}}
\newcommand{\ed}{\end{displaymath}}
\newcommand{\bp}{\underline{\bf Proof}:\ }
\newcommand{\ep}{{\hfill $\Box$}\\ }
\newtheorem{1}{LEMMA}[section]
\newtheorem{2}{THEOREM}[section]
\newtheorem{3}{COROLLARY}[section]
\newtheorem{4}{PROPOSITION}[section]
\newtheorem{5}{REMARK}[section]
\newtheorem{20}{OBSERVATION}[section]
\newtheorem{10}{DEFINITION}[section]
\newtheorem{30}{RESULTS}[section]
\newtheorem{40}{CLAIM}[section]
\newtheorem{50}{ASSUMPTION}[section]
\newtheorem{60}{EXAMPLE}[section]
\newtheorem{70}{ALGORITHM}[section]
\newtheorem{80}{PROBLEM}
\newcommand{\be}{\begin{equation}}
\newcommand{\ee}{\end{equation}}
\newcommand{\ba}{\begin{array}}
\newcommand{\ea}{\end{array}}
\newcommand{\bea}{\begin{eqnarray}}
\newcommand{\eea}{\end{eqnarray}}
\newcommand{\bqn}{\begin{eqnarray*}}
\newcommand{\eqn}{\end{eqnarray*}}


\newcommand{\e} { \ = \ }
\newcommand{\leqs}{ \ \leq \ }
\newcommand{\geqs}{ \ \geq \ }
\def\theequation{\thesection.\arabic{equation}}
\def\bReff#1{{\bRm
(\bRef{#1})}}
\newcommand{\eps}{\varepsilon}
\newcommand{\sgn}{\operatorname{sgn}}
\newcommand{\sign}{\operatorname{sign}}
\newcommand{\Vol}{\operatorname{Vol}}
\newcommand{\Var}{\operatorname{Var}}
\newcommand{\Cov}{\operatorname{Cov}}
\newcommand{\vol}{\operatorname{vol}}
\newcommand{\var}{\operatorname{var}}
\newcommand{\cov}{\operatorname{cov}}
\renewcommand{\Re}{\operatorname{Re}}
\renewcommand{\Im}{\operatorname{Im}}
\newcommand{\bE}{{\mathbb E}}
\newcommand{\bR}{\mathbb{R}}
\newcommand{\bN}{{\mathbb N}}
\newcommand{\bC}{\mathbb{C}}
\newcommand{\bF}{\mathbb{F}}
\newcommand{\bQ}{{\mathbb Q}}
\newcommand{\bZ}{{\mathbb Z}}
\newcommand{\cA}{{\mathcal A}}
\newcommand{\cB}{{\mathcal B}}
\newcommand{\cC}{{\mathcal C}}
\newcommand{\cD}{{\mathcal D}}
\newcommand{\cE}{{\mathcal E}}
\newcommand{\cF}{{\mathcal F}}
\newcommand{\cG}{{\mathcal G}}
\newcommand{\cH}{{\mathcal H}}
\newcommand{\cI}{{\mathcal I}}
\newcommand{\cJ}{{\mathcal J}}
\newcommand{\cK}{{\mathcal K}}
\newcommand{\cL}{{\mathcal L}}
\newcommand{\cM}{{\mathcal M}}
\newcommand{\cN}{{\mathcal N}}
\newcommand{\cO}{{\mathcal O}}
\newcommand{\cP}{{\mathcal P}}
\newcommand{\cQ}{{\mathcal Q}}
\newcommand{\cR}{{\mathcal R}}
\newcommand{\cS}{{\mathcal S}}
\newcommand{\cT}{{\mathcal T}}
\newcommand{\cU}{{\mathcal U}}
\newcommand{\cV}{{\mathcal V}}
\newcommand{\cW}{{\mathcal W}}
\newcommand{\cX}{{\mathcal X}}
\newcommand{\cY}{{\mathcal Y}}
\newcommand{\cZ}{{\mathcal Z}}
\newcommand{\bx}{{\mathbf x}}
\newcommand{\by}{{\mathbf y}}
\newcommand{\bz}{{\mathbf z}}
\newcommand{\bba}{{\mathbf a}}
\newcommand{\bbb}{{\mathbf b}}
\newcommand{\bbc}{{\mathbf c}}


\title{A homotopy method for solving multilinear systems with M-tensors}
\author{ 
Lixing Han \thanks{ Department of Mathematics,  University of Michigan-Flint, Flint, MI 48502, USA.
Email: \texttt{lxhan@umflint.edu}. The author was supported in part by a 60th Anniversary Research Grant, Office of the Provost, UM-Flint.} 
}
\date{December 22, 2016}
\maketitle

\begin{abstract}  
 Multilinear systems of equations arise in various applications, such as numerical partial differential equations, data mining, and tensor complementarity problems. In this paper,  we propose a homotopy method for finding the unique positive solution to a multilinear system with a nonsingular M-tensor and a positive right side vector.  We analyze the method and prove its convergence to the desired solution. We report some numerical results based on an implementation of the proposed method using a prediction-correction approach for path following.
 
\end{abstract}



\section{Introduction}
\label{Intro}
\setcounter{equation}{0}
Let $\bR$ and $\bC$ be the real field and complex field, respectively. We denote the set of all $m$th-order, $n$-dimensional real tensors  by $\bR^{[m,n]}$.  
For a tensor $\cA \in \bR^{[m,n]}$ and a vector $b \in \bR^n$, a multilinear system is defined as 
\be
\label{multisys}
\cA x^{m-1} = b,
\ee
where $x=[x_1, x_2, \ldots, x_n]^T $ is the unknown vector, and  $\cA x^{m-1}$ denotes the column vector  whose $i$th entry is
$$
(\cA x^{m-1})_i = \sum_{i_2, \cdots, i_m =1}^{n} A_{i  i_{2}  \cdots  i_m}x_{i_2} \cdots x_{i_m},
$$ 
for $i=1,2, \ldots, n$.  Multilinear systems of the form (\ref{multisys}) arise in a number of applications, such as numerical partial differential equations, data mining, and tensor complementarity problems (see for example, \cite{DW16, LN15, LQX15}).

 In their pioneering works, Qi \cite{Qi05} and Lim \cite{Lim05} independently introduced the concept of tensor eigenvalues. We say that $(\lambda, x) \in \bC \times \bC^{n} \backslash \{0\}$ is an eigenpair of a tensor $\cA \in \bR^{[m,n]}$ if 
$$
\cA x^{m-1} = \lambda x^{[m-1]},
$$
where $x^{[m-1]} = [x_{1}^{m-1}, x_{2}^{m-1}, \ldots, x_{n}^{m-1}]^T$.    Let $\rho (\cA)$  denote the spectral radius of tensor $\cA$, that is, 
$$
\rho (\cA)  = \max \{ |\lambda|: \lambda {\rm \ is \ an \ eigenvalue \ of} \cA. \}.
$$

 Recently, the notion of M-tensors has been proposed and their properties have been studied in \cite{DQW13, ZQZ14}. A tensor $\cA \in \bR^{[m,n]}$ is called an M-tensor if it can be written as $\cA = s \cI -B$, in which $ \cI$ is the $m$th-order, $n$-dimensional identity tensor, $\cB$ is a nonnegative tensor (that is, each entry of $\cB$ is nonnegative), and $s \geq \rho (\cB)$. Furthermore, $\cA$ is called a nonsingular M-tensor if $s > \rho(\cB)$.    

In \cite{DW16}, Ding and Wei investigated the solutions of the multilinear system (\ref{multisys}) when  the coefficient tensor $\cA$ is an M-tensor. In particular, they show that the system (\ref{multisys}) has a unique positive solution if $\cA$ a nonsingular M-tensor and $b$ is a positive vector (see \cite[Theorem 3.2]{DW16}).  They generalized the Jacobi and Gauss-Seidel methods for linear systems to find the unique positive solution of the multilinear system (\ref{multisys}). They also proposed to use Newton's method when the nonsingular M-tensor $\cA$ is symmetric and numerically showed that Newton's method is much faster than the other methods. However,  it is unclear whether or not the Newton method proposed in \cite{DW16} always works when $\cA$ is not symmetric.   

In this paper, we propose a homotopy method for finding the unique positive solution of the multilinear system (\ref{multisys}) and prove its convergence. The homotopy method is implemented using an Euler-Newton prediction-correction approach for path tracking. Numerical experiments show the efficiency of our method.

The paper is organized as follows.  We  introduce our homotopy method and prove its convergence in Section 2. Then we give some numerical results in Section 3.

\section{A Homotopy Method}
\label{sec2}
\setcounter{equation}{0}

We are to design a homotopy method for  finding the unique positive solution of the system (\ref{multisys}) when $\cA$ is a nonsingular M-tensor and $b$ is a positive vector.  For this purpose, we will solve the following polynomial system:
\be
\label{target}
P(x) =  \cA x^{m-1} - b = 0.
\ee
We choose the starting system
\be
\label{start}
Q(x) = \cI x^{m-1} - b = 0,
\ee
and construct the following homotopy 
\be
\label{homotopy}
H(x,t) = (1-t) Q(x) + t P(x) =0,  \ \ \   t \in [0,1].
\ee
Note that the starting system (\ref{start}) trivially has a unique positive solution 
\be
 x_0=[b_{1}^{1/(m-1)}, b_{2}^{1/(m-1)}, \ldots,b_{n}^{1/(m-1)}]^T.
\ee
Moreover,  the homotopy $H(x,t)$ can be expressed as 
\be
\label{hexpress}
H(x,t) = (t \cA +(1-t) \cI) x^{m-1} - b.
\ee

The partial derivatives matrix $D_x H(x,t)$ of the homotopy $H(x, t)$ plays an important role in our method. To compute this matrix, we need to partially symmetrize  tensor $\cA=(A_{i_1,i_2,\ldots, i_m})$ with respect to the indices $i_2, \ldots, i_m$. Specifically, we define the partially symmetrized tensor 
$\hat{\cA}=(\hat{A}_{i_1, i_2, \ldots, i_m})$ by
\be
\label{psymm}
\hat{A}_{i_1 i_2 \ldots i_m} = \frac{1}{(m-1)!} \sum_{\pi} A_{i_1 \pi(i_2\ldots   i_m)},
\ee
 where the sum is over all the permutations $\pi(i_2\ldots  i_m)$.  The following lemma shows that this partial symmetrization preserves the nonsingular M-tensor structure. 

\begin{1} 
\label{lemma1}
If $\cA \in \cR^{[m,n]}$ is a nonsingular M-tensor, so is $\hat{\cA}$.  
\end{1}
\bp
Since $\cA$ is an M-tensor, there is a nonnegative tensor $\cB$ such that $\cA= s \cI - \cB$. Then $\hat{\cB}$ is nonnegative and moreover, $\hat{\cA} = s \cI - \hat{\cB}$.
According to \cite[Page 3277, Conditions D1 and D4]{DQW13}, $\cA$ is a nonsingular M-tensor if and only if there is a positive vector $y \in \bR$ such that $\cA y^{m-1}$ is a positive vector.  Note that 
$\cA x^{m-1} = \hat{\cA} x^{m-1}$ for all $x \in \bR^n$. Thus, $\hat{\cA} y^{m-1}$ is a positive vector. 
  Using the results in \cite{DQW13} again, $\hat{\cA}$ is a nonsingular M-tensor.  
\ep

The partial derivatives matrix of $\cA x^{m-1}$ with respect to $x$ is 
\be
\label{jacobian}
D_x \cA x^{m-1} = (m-1) \hat{\cA} x^{m-2}. 
\ee
 Therefore, the partial derivatives of $H$ with respect to  $x$ and $t$ are:
$$
D_{x} H (x,t)  =
(m-1) ( t \hat{\cA} +(1-t) \cI) x^{m-2},
$$
and 
$$
D_{t} H (x,t)  = 
 (\cA-\cI)x^{m-1},
$$
respectively.

\begin{2}
\label{thm1} 
Suppose that $\cA   \in \bR^{[m,n]}$ is a nonsingular M-tensor and $b \in \bR^n$ is a positive vector.  Then there exists a number $\tau_0 >0$ such that, for each $t \in [0,1+\tau_0)$, \\
(a)  $H(x, t)=0$ has a unique positive solution $ x(t)$; \\
(b)  the partial derivatives matrix 
$$
D_{x} H(x(t), t) = (m-1) ( t \hat{\cA} +(1-t) \cI ) x(t)^{m-2}
$$ 
is nonsingular. 
\end{2}
\bp
Let $\cA = s \cI - \cB$ such that $\cB$ is a nonnegative tensor and $s > \rho(\cB)$. Note that
 $$
t \cA +(1-t) \cI =  (st+1-t)\cI -t\cB.
$$
Clearly, this tensor is a nonsingular M-tensor for $ 0\leq t \leq 1$. Choose  
\[
\tau_0 = \left \{  \begin{array}{l}
\frac{s-\rho(\cB)}{\rho(\cB)-s+2 }, \ \  {\rm if} \ \ \rho(\cB)-s+2 >0, \\
1, \hskip 0.6in  {\rm if} \ \ \rho(\cB)-s+2  \leq 0.
\end{array}
\right.
\] 
Then $\tau_0>0$ and $ \displaystyle{\frac{st+1-t}{t}  > \frac{s+\rho(\cB)}{2} > \rho(\cB)}$  for $ 1 \leq t < 1+\tau_0$. This   implies that  
$$\displaystyle{t \cA +(1-t) \cI = t \left ( \frac{st+1-t}{t} \cI - \cB \right )}
$$
 is a nonsingular M-tensor when $ 1 \leq t < 1+\tau_0$. 
Therefore,  $t \cA +(1-t) \cI $ is a nonsingular M-tensor for each $ t \in [0, 1+\tau_0)$. It follows that
 $H(x,t) = (t \cA +(1-t) \cI) x^{m-1} - b =0$
has a unique positive solution $x(t)$ for each $ t \in [0, 1+\tau_0)$ 
by \cite[Theorem 3.2]{DW16}.  Moreover, $t \hat{\cA} +(1-t) \cI $ is a nonsingular M-tensor for each $ t \in [0, 1+\tau_0)$ by Lemma \ref{lemma1}. Note that the matrix
  $ (t \hat{\cA} +(1-t) \cI) x(t)^{m-2}$ is Z-matrix, and 
$$
[ (t \hat{\cA} +(1-t) \cI) x(t)^{m-2} ] x(t)= (t \hat{\cA} +(1-t) \cI) x(t)^{m-1} = (t \cA +(1-t) \cI) x(t)^{m-1} = b 
$$  
is a positive vector, we must have that $D_xH(x(t),t)=(m-1)(t \hat{\cA} +(1-t) \cI) x(t)^{m-2}$ is nonsingular by \cite[Chapter 6]{BP94}.
\ep

\begin{3}
\label{cor1}
Suppose that $\cA   \in \bR^{[m,n]}$ is a nonsingular M-tensor and $b \in \bR^n$ is a positive vector. Then the positive solutions $x(t)$ of $H(x,t)=0$ for
$t \in [0,1+\tau_0)$  form a smooth curve in $\bR_{++}^n$, where  $\bR_{++}^n$ is the set of positive $n$-vectors. 
\end{3}
\bp
By Theorem \ref{thm1}, $H(x,t)=0$ has a unique positive solution $x(t)$ for each $t \in [0,1+\tau_0)$.   The conclusion follows by using  the Implicit Function Theorem and a continuation argument  (\cite{Li15}). 
\ep

Our next theorem shows that the homotopy (\ref{homotopy}) works.

\begin{2}
\label{main}
Suppose that $\cA \in \bR^{[m,n]}$ is a nonsingular M-tensor and $b \in \bR^n$ is a positive vector.  Starting from the initial 
$x(0) = x_0 = [b_{1}^{1/(m-1)}, b_{2}^{1/(m-1)}, \ldots,b_{n}^{1/(m-1)}]^T$,  let $ x(t)$ be the solution curve obtained by solving the homotopy $H(x,t)=0$ in $\bR_{++}^{n} \times [0,1]$. Then  $ x(1)$ is the unique positive solution 
of the system (\ref{multisys}). 
\end{2}
\bp
According to Corollary \ref{cor1}, the positive solutions $x(t)$ of $ H( x, t)=0 $ for $t \in [0, 1+\tau_0)$ form a smooth curve in  $\bR_{++}^{n}$. 
Differentiating $H(x(t), t)=0$ with respect to $t$ gives
\be
\label{diffeq}
D_{x}  H(x(t),t)    \cdot  \frac{d x}{d t}
= - D_t H(x(t),t).
\ee
Since $D_{x}  H (x(t),t)$ is nonsingular for all $t \in [0,1+\tau_0)$, this system of differential equations is well defined for $t  \in [0,1+\tau_0)$. 
We can follow the curve by solving this system with the initial condition $x(0)=x_0$. Clearly,  $ x(1)$ is the unique positive solution 
of the system (\ref{multisys}). 
\ep

We now present our homotopy method for finding the unique positive solution of (\ref{multisys}). 

\begin{framed}

\begin{70} {Finding the positive solution for (\ref{multisys}) when $\cA$ is a nonsingular M-tensor and $b$ a positive vector.}
\label{algorithmA}

\ \\
{\bf  \large Initialization.}  Choose initial $x_0= [b_{1}^{1/(m-1)}, b_{2}^{1/(m-1)}, \ldots,b_{n}^{1/(m-1)}]^T$.  \\
\ \\
{\bf  \large Path following.}  Solve the differential system (\ref{diffeq}) with the initial condition $x(0) = x_0$ in $\bR_{++}^{n}$.  $x(1)$ is the desired solution for system (\ref{multisys}).

\end{70}
\end{framed}

\section{Numerical Results}
\label{sec3}
\setcounter{equation}{0}

We have implemented Algorithm \ref{algorithmA} in Matlab.  The code can be downloaded from: 
\begin{verbatim}
    http://homepages.umflint.edu/~lxhan/software.html
\end{verbatim} 

In our implementation,  an Euler-Newton type predication-correction approach (\cite{AG90, SW05}) with an adaptive stepsize for solving the system of differential equations (\ref{diffeq}) with  initial condition $x(0)=x_0$ is used. Moreover, it solves the scaled system 
$$
\bar{\cA} x^{m-1} = \bar{b},
$$
where $\bar{\cA} = \cA / \omega$, $\bar{b} = b/\omega$, and $\omega$ is the largest value among the entries of $b$ and the absolute values of entries of $\cA$.  The code terminates if the residue of the scaled system 
$$
\|\bar{\cA} x(1)^{m-1} - \bar{b} \|_2 \leq 10^{-12}.
$$

To test the effectiveness of Algorithm \ref{algorithmA}, we did some numerical experiments.  All the experiments were done using MATLAB 2014b on a laptop computer with Intel Core i7-4600U at 2.10 GHz and 8 GB memory running Microsoft Windows 7.  The tensor toolbox of \cite{BK15} was used to compute tensor-vector products and to compute partially symmetrized tensor $\hat{\cA}$.  

The examples we tested are generated by the following method. We chose nonnegative tensor $\cB \in \bR^{[m,n]}$ whose entires are  uniformly distributed in $(0,1)$. Set 
$$
s= (1+\epsilon) \cdot \max_{1\leq i \leq n} \sum_{i_2, \ldots, i_m} \cB_{i,i_2,\ldots,i_m},
$$
for some $\epsilon>0$.  Let $\cA= s \cI - \cB$. According to \cite{Qi05, CQZ13}, 
$$
\rho(\cB) \leq \max_{1\leq i \leq n} \sum_{i_2, \ldots, i_m} \cB_{i,i_2,\ldots,i_m}.
$$
Thus, $s> \rho(\cB)$ and $\cA$ is a nonsingular M-tensor.  We chose the right side vector $b \in \bR^n$ with entires  uniformly distributed in $(0,1)$.

\begin{table}[htbp]
\label{table1}
\begin{center}
\begin{tabular}{|c|c|c|c|c|}
\hline
$(m,n)$  & euitr & nwitr & time  & residue \\ \hline
(3,10) & 5  & 11 &  0.098  &  $ 3.8531 \times 10^{-15}   $  \\ \hline
(3,50) & 5  & 10 &  0.126  &  $ 3.9837 \times 10^{-14}   $  \\ \hline
(3,100) & 5  & 9 &  0.289  &  $ 5.8008 \times 10^{-13}  $  \\ \hline
(3,200) & 5  & 8 &  0.929  &  $ 2.9054 \times 10^{-11}  $  \\ \hline
(3,400) & 5  & 7 &  8.099  &  $ 1.3199 \times 10^{-8}  $  \\ \hline
(4,10) & 5  & 10 &  0.134  &  $ 2.2176 \times 10^{-12}  $  \\ \hline
(4,50) & 5  & 8 &  1.019  &  $ 1.6437 \times 10^{-8}  $  \\ \hline
(4,80) & 5  & 8 &  8.902  &  $ 3.3791 \times 10^{-11}  $  \\ \hline
(4,100) & 5  & 8 &  19.423  &  $ 1.3962 \times 10^{-10}  $  \\ \hline
(5,10) & 5  & 9 &  0.165  &  $ 7.0535 \times 10^{-13}  $  \\ \hline
(5,20) & 5  & 10 &  1.646  &  $ 2.4450 \times 10^{-13}  $  \\ \hline
(5,40) & 5  & 9 &  55.656  &  $ 5.4560 \times 10^{-10}  $  \\ \hline
(6,5)  & 5  & 10 &  0.242 &  $ 1.3524 \times 10^{-12}  $  \\ \hline
(6,10) & 5  & 9 &  1.483 &  $ 3.1020 \times 10^{-8}  $  \\ \hline
(6,15) & 5  & 10 &  31.232  &  $ 6.2806 \times 10^{-11}  $  \\ \hline
\end{tabular}
\caption{Numerical Results} 
\end{center}
\end{table}

We tested the algorithm on tensors of various sizes by choosing  different values of $m$ and $n$. We used $\epsilon = 0.01$ as in \cite{DW16}. We now summarize the numerical results in Table 1. In this table,  \texttt{euitr} and \texttt{nwitr} denote the number of Euler prediction steps and the total number of Newton correction steps were used, \texttt{time} denotes the CPU time used (in seconds) when the algorithm terminated, and  \texttt{residue} denotes the residue $\|\cA x(1)^{m-1}-b\|_2$ of the original system (\ref{multisys}) at termination.  

For the examples we tested, we observe that Algorithm \ref{algorithmA} can find the positive solution of the multilinear system (\ref{multisys}) when $\cA$ is a nonsingular M-tensor and $b$ is a positive vector. It is efficient in terms of both \texttt{euitr} and \texttt{nwitr}.  We remark that for the $(4,100), (5,40), (60, 15)$ cases, the relatively large CPU time used by Algorithm \ref{algorithmA} is mainly due to the procedure of partially symmetrizing tensor $\cA$.  A more efficient symmetrization method can help save the CPU time for such cases. Of course,  the partial symmetrization is not needed if the tensor $\cA$ is symmetric. Overall, the numerical results show that Algorithm \ref{algorithmA} is quite promising.


\ \\



\begin{thebibliography}{999}

\bibitem{AG90} E.L. Allgower and K. Georg, {\it Numerical Continuation Methods, an Introduction}, Springer Series in Comput. Math., Vol 13, Springer-Verlag (Berlin, Heidelberg, New York), 1990.

\bibitem{BK15} B.W. Bader, T.G. Kolda and others, MATLAB Tensor Toolbox Version 2.6, 2015.   

\bibitem{BP94} A. Berman and R.J. Plemmons, {\it Nonnegative Matrices in the Mathematical Sciences}, Classics in Applied Mathematics, SIAM, Philadelphia, 1994. 

\bibitem{CQZ13} K.C. Chang, L. Qi, and T. Zhang, A survey of the spectral theory of nonnegative tensors,  {\it Numerical Linear Algebra with Applications}, 2013, 20: 891--912.

\bibitem{DQW13} W. Ding, L. Qi, and Y. Wei, M-tensors and nonsingular M-tensors, {\it Linear Algebra and its Applications}, 2013, 439(10): 3264--3278.

\bibitem{DW16} W. Ding and Y. Wei, Solving multi-linear systems with M-tensors, {\it Journal of Scientific Computing}, 2016, 68: 689--715. 

\bibitem{LN15} X. Li and M.K. Ng, Solving sparse non-negative tensor equations: algorithms and applications, {\it Frontiers of Mathematics in China}, 2015, 10(3): 649--680.

\bibitem {Li15} T.Y. Li, Homotopy methods,  in {\it Encyclopedia of Applied and Computational Mathematics}, B. Engquist, ed., Springer, Berlin, 2015: 653--656.

\bibitem{Lim05} L.-H. Lim, Singular values and eigenvalues of tensors: a variational approach, {\it Proceedings of the IEEE International Workshop on Computational Advances
in Multi-Sensor Adaptive Processing (CAMSAP'05)}, 2005, 1: 129--132.

\bibitem{LQX15} Z. Luo, L. Qi, and N. Xiu, The sparsest solutions to Z-tensor complementarity problems, arXiv: 1505.00993, 2015.  

\bibitem{Qi05} L. Qi, Eigenvalues of a real supersymmetric tensor, {\it Journal of Symbolic Computation}, 2005, 40:
1302--1324.

\bibitem{SW05} A.J. Sommese and W.W. Wampler, {\it The Numerical Solution of Systems of Polynomials Arising
in Engineering And Science} World Scientific Pub Co Inc, 2005.

\bibitem{ZQZ14} L. Zhang, L. Qi, and G. Zhou, M-tensors and some applications, {\it SIAM Journal on Matrix Analysis and Applications}, 2014, 35(2): 437--452. 

\end{thebibliography}
\end{document}